\documentclass{llncs} 
\usepackage[colorlinks=true,
linkcolor=webgreen,
filecolor=webbrown,
citecolor=webgreen]{hyperref}

\definecolor{webgreen}{rgb}{0,.5,0}
\definecolor{webbrown}{rgb}{.6,0,0}

\usepackage{psfig,amssymb}
\def\binom#1#2{{#1}\choose{#2}}
\newcommand{\eqn}[1]{(\ref{#1})}

\begin{document}
\title{My Favorite Integer Sequences}
\author{N. J. A. Sloane}
\institute{Information Sciences Research \\
           AT\&T Shannon Lab \\
           Florham Park, NJ 07932-0971 USA \\
           Email: \href{mailto:njas@research.att.com}{njas@research.att.com}
}
\date{}
\maketitle 

\begin{abstract}
This paper gives a brief description of the author's database of integer sequences,
now over 35 years old, together with a selection of a few of the most interesting
sequences in the table.
Many unsolved problems are mentioned.

\vspace{.5\baselineskip}
This paper was published (in a somewhat different form) in
{\em Sequences and their Applications (Proceedings of SETA '98)}, C. Ding, T. Helleseth and H. Niederreiter (editors), Springer-Verlag, London, 1999, pp. 103-130. 
Enhanced pdf version prepared Apr 28, 2000.

\vspace{.5\baselineskip}
The paragraph on ``Sorting by prefix reversal'' on page \pageref{sorting}
was revised Jan. 17, 2001.
\end{abstract}

\section{How it all began}
I started collecting integer sequences in December 1963, when I was a graduate student at
Cornell University, working on perceptrons (or what are now called neural networks).
Many graph-theoretic questions had arisen, one of the simplest
of which was the following.

Choose one of the $n^{n-1}$ rooted labeled trees with $n$ nodes at random, and pick a random node:
what is its expected height above the root?
To get an integer sequence, let $a_n$ be the sum of the heights of all nodes in all
trees, and let $W_n = a_n/n$.
The first few values $W_1, W_2 , \ldots$ are
\begin{center}
0, 1, 8, 78, 944, 13800, 237432, $\ldots$,
\end{center}
a \htmladdnormallink{sequence}{http://www.research.att.com/cgi-bin/access.cgi/as/njas/sequences/eisA.cgi?Anum=000435} engraved on my memory.
I was able to calculate about ten terms, but I needed to know how
$W_n$ grew in comparison with $n^n$, and it was impossible to guess this from 
so few terms.
So instead I tried to guess a formula for the $n$th term.
Again I was unsuccessful.
Nor could I find this sequence in Riordan's book \cite{R1}, although
there were many sequences that somewhat resembled it.

So I started collecting all the sequences I could find, entering them on punched cards,
thinking that if any of these sequences came up in another problem,
at least I would know what {\em they} were.

I never did find that sequence in the literature, but I learned P\'{o}lya's theory of counting
and (with John Riordan's help) obtained the answer, which appears in \cite{RS1}.
There {\em is} a simple formula, although maybe not simple enough to be guessed:
$$W_n = (n-1)! \sum_{k=0}^{n-2} \frac{n^k}{k!} ~.$$
An old formula of Ramanujan \cite{Wat29} then implies that
$$\frac{W_n}{n^n} \sim \sqrt{\frac{2\pi}{n}} , \quad\mbox{as}\quad n \to \infty~,
$$
which was what I needed.
That sequence became number \htmladdnormallink{A435}{http://www.research.att.com/cgi-bin/access.cgi/as/njas/sequences/eisA.cgi?Anum=000435} in the collection.

The idea of a ``dictionary'' of integer sequences was received with enthusiasm by many
people, and in 1973 I published \cite{HIS}, containing about 2400 sequences,
arranged lexicographically.
One correspondent commented on the book by saying
``There's the Old Testament, the New Testament and the Handbook of Integer Sequences''.

Over the next twenty years an enormous amount of new material arrived,
preprints, reprints, postcards, typewritten letters, handwritten letters, etc., and it was
not until 1995 that --- with Simon Plouffe's help ---
a sequel \cite{EIS} appeared.
This contained 5500 sequences.

Around the same time I set up two services that can be used to consult the database
via electronic mail.
The first has the address
\href{mailto:sequences@research.att.com}{sequences@research.att.com},
and simply looks up a sequence in the table.
The second
email address, which is
\href{mailto:superseeker@research.att.com}{superseeker@research.att.com},
tries hard to find an explanation even for a sequence not in the table.

A large number of new sequences started arriving as soon as \cite{EIS}
appeared, and when in 1996 the total number reached 16,000, three times the number in the book,
I decided to set up a web site for the database \cite{OEIS}.
New sequences still continue to pour in, at about 10000 per year.
At present, in April 1999, the database contains about 50000 sequences, and the web site
receives over 2500 hits per day.

The main reason for this rapid expansion is that in the two books I only included sequences
that had been or were about to be published.
For the on-line version, where storage space is no longer a limitation, any well-defined
and sufficiently interesting sequence is eligible for inclusion.

There is now also an electronic {\em Journal of Integer Sequences}
\cite{JIS} and a mailing list for sequence fans \cite{Fan}.

The following sections describe how the database is used (Section 2) and the kinds of sequences
it contains (Section 3).
Section 4 discusses a few ``hard'' sequences and Section 5 some recursive examples.
Then Sections 6--8 describe sequences associated with meandering rivers
and stamp-folding, extremal codes and lattices, and Levine's sequence.

Besides integer sequences, the table also contains many examples of arrays of numbers,
Pascal's triangle, for example --- see Section 9.
The final three sections discuss sequences associated with the Wythoff array, the boustrophedon transformation of sequences,
and Tchoukaillon solitaire.

\section{How the database is used}
The main applications of the database are in identifying sequences or in finding out the
current status of a known sequence.

The database has been called a mathematical analogue of a ``fingerprint file''
\cite{Cipra}.
You encounter a number sequence, and wish to know if anyone has ever come across it before.
If your sequence is in the database, the reply will provide a description, the first 50
or so terms (usually enough to fill three lines on the screen), and,
when available, formulae, recurrences,
generating functions, references,
computer code for producing the sequence, links to relevant web sites, etc.

Let me illustrate how the database is used with a typical story.
Last summer the following question arose at AT\&T Labs in connection with a quantization
problem \cite{SVS}.
Given an $n$-dimensional lattice $\Lambda$, for which integers $N$ does $\Lambda$ have a sublattice of index $N$ that is geometrically similar to $N$?

For the two-dimensional root lattice $A_2$, for example, it is easy to see that a necessary
and sufficient condition is that $N$ be of the form $a^2 + ab + b^2$.

However, for the four-dimensional lattice $A_4$ the situation is more complicated.
The first thing we did was to run a computer search, which showed that $A_4$ has a similar
sublattice of index $N= c^2$ if and only if $c$ is one of the numbers
$$1,4,5,9,11,16,19,20,25,29,31,36, \ldots$$
This turned out to be sequence \htmladdnormallink{A31363}{http://www.research.att.com/cgi-bin/access.cgi/as/njas/sequences/eisA.cgi?Anum=031363} in the database, with a reference to Baake \cite{Baa1},
where it appears
in an apparently different context
as the indices of coincidence site sublattices in a certain
three-dimensional quasicrystal.
\cite{Baa1} identifies these numbers as those positive integers in which all primes congruent to 2 or 3 (modulo 5) appear to an even power.
This was a very useful hint in getting started on our problem.
We were able to show that this is the correct condition for the $A_4$ lattice, and to find analogous results for a number of
other lattices \cite{CRS}.
(However, we have not yet found a direct connection between $A_4$ and the quasicrystal problem.
Nevertheless, the occurrence of the same numbers in the two problems cannot be entirely
coincidental.)

My files contain many similar stories.
Some other examples can be found in Chapter 3 of \cite{EIS}.

Most of the applications however are less dramatic.
One encounters a sequence in the middle of a calculation,
\htmladdnormallink{perhaps}{http://www.research.att.com/cgi-bin/access.cgi/as/njas/sequences/eisA.cgi?Anum=001405}
\begin{center}
1,1,2,3,6,10,20,35,70,126, 252, $\ldots$,
\end{center}
and one wants to know quickly what it is --- preferably a formula (in this case it is ${\binom{n}{[n/2]}}$)
or generating function.

It is worth emphasizing a special case of this: the
simplification of binomial coefficient sums.
Powerful methods are available for simplifying such sums by
computer \cite{Nem}, \cite{AB}.
But if one is in a hurry, one can first try evaluating the initial terms and looking up the sequence
in the table.
E.g. the sum
$$a(n) = \sum_{k=0}^n {\binom{2n-2k}{n-k}}^2 {\binom{2k}{k}}$$
produces
$$1, 8,88,1088, 14296, 195008 , \ldots$$
(\htmladdnormallink{A36917}{http://www.research.att.com/cgi-bin/access.cgi/as/njas/sequences/eisA.cgi?Anum=036917}), and the entry in the table supplies a recurrence
$$n^3 a(n) = 16 \left( n - \frac{1}{2} \right)
(2n^2 -2n+1) a(n-1) - 256 (n-1)^3 a(n-2)$$
and a reference (\cite{AB}, p. ix, although the sum is misstated in the first printing).
I have begun entering
into the database
all the sequences corresponding to the singly-indexed identities in Gould's table \cite{Gould}
of binomial coefficient identities.

A related application is in identifying arithmetic inequalities.
Suppose you suspect that
\begin{equation}\label{EqIa}
\sigma (n) \ge d(n) + \phi (n) , \quad\mbox{for}\quad
n \ge 2 ~,
\end{equation}
where $\sigma, d$ and $\phi$ are respectively the sum of divisors, number of divisors, and Euler totient functions.
You could evaluate the sequence formed by the difference of the two sides, which for $n \ge 1$ is
$$-1, 0,0,2,0,6,0,7,4,10,0,18,0,14,12, \ldots$$
(\htmladdnormallink{A46520}{http://www.research.att.com/cgi-bin/access.cgi/as/njas/sequences/eisA.cgi?Anum=046520}).
The table then points you to a reference (\cite{HNT}, \S I.3.1) where this is stated as a theorem.
(Again I would like to get more examples of such sequences.)

Another important use for the database is in finding out the present state of knowledge about
some problem.
A second story will illustrate this.
The number of Latin squares of order $n$ $(Q_n )$ is given by sequence \htmladdnormallink{A315}{http://www.research.att.com/cgi-bin/access.cgi/as/njas/sequences/eisA.cgi?Anum=000315}
(see Fig.~\ref{FLS1}).
Computing $Q_n$ is one of the famous hard problems in combinatorics.
\begin{figure}[htb]
$$
\begin{array}{r@{~~~~}r}
n & \multicolumn{1}{c}{Q_n} \\ \hline
1 & 1 \\
2 & 1 \\
3 & 1 \\
4 & 4 \\
5 & 56 \\
6 & 9408 \\
7 & 16942080 \\
8 & 535281401856 \\
9 & 377597570964258816 \\
10 & 7580721483160132811489280
\end{array}
$$

\caption{Sequence \protect\htmladdnormallink{A315}{http://www.research.att.com/cgi-bin/access.cgi/as/njas/sequences/eisA.cgi?Anum=000315}, the number of Latin squares of order $n$ (McKay and Rogoyski \cite{MR95}).}
\label{FLS1}
\end{figure}

In 1991 Brendan McKay, at the Australian National University in Camberra,
computed $Q_{10}$ (which as can be seen is a rather large number).
When he checked the database he found that the same value had recently been obtained by
Eric Rogoyski of Cadence Design Systems in San Jose, California.
As it turned out the two methods were similar but not identical, and he and Rogoyski ended up writing a joint paper
\cite{MR95}, acknowledging the database for bringing them together.

\section{Types of Sequences}
The database contains sequences from all branches of science, including
\begin{itemize}
\item
enumeration problems (combinatorics, graph theory, lattices, etc.)
\item
number theory (number of solutions to $x^2 + y^2 + z^2 =n$, etc.)
\item
game theory (winning positions, etc.)
\item
physics (paths on lattices, etc.)
\item
chemistry (sizes of clusters of atoms, etc.)
\item
computer science (number of steps to sort $n$ things, etc.)
\item
communications ($m$-sequences, weight distributions of codes, etc.)
\item
puzzles
\item
etc.
\end{itemize}

To be accepted into the database, a sequence must be well-defined and interesting.
However, when in doubt, my tendency is to accept rather than to reject.
The amazing coincidences of the Monstrous Moonshine investigations \cite{MM}
make it difficult to say that a particular sequence, no matter how obscure,
will never be of interest.

Readers are urged to send me any sequences they come across that are not at present in the
database.
There is a convenient electronic form for this purpose in \cite{OEIS}.
Some of the reasons for sending in your sequence are as follows:
\begin{itemize}
\item
this stakes your claim to it
\item
your name is immortalized
\item
the next person who comes across it will be grateful and, not least,
\item
you may benefit from this yourself, when you come across the same sequence some weeks from now.
\end{itemize}

Often one finds that a particular project may involve dozens of sequences, all variants of a few
basic ones.
Ideally you should send them all to the database!

\section{Hard sequences}
One of the keywords used in the database is ``hard'',
which indicates that the term following those given is not known.
Besides the Latin squares problem mentioned above, some other classic hard sequences are the following: \\

\noindent{\bf Projective planes.}
The number of projective planes of orders $n=2,3, \ldots, 10$ (\htmladdnormallink{A1231}{http://www.research.att.com/cgi-bin/access.cgi/as/njas/sequences/eisA.cgi?Anum=001231}):
$$1,1,1,1,0,1,1,4,0 ,$$
where the last term refers to the result of Clement Lam et~al. \cite{LamTS89}
(completing work begun in \cite{MacST73})
that there is no projective plane of order 10. \\

\noindent{\bf The Poincar\'{e} conjecture.}
The number of differential structures on the $n$-sphere, for
$n=1,2, \ldots, 16$ (\htmladdnormallink{A1676}{http://www.research.att.com/cgi-bin/access.cgi/as/njas/sequences/eisA.cgi?Anum=001676}):
$$1,1,1?, 1,1,1, 28, 2, 8, 6, 992, 1,3,2, 16256,2,$$
as given by Kervaire and Milnor \cite{KM63}.
The {\em Poincar\'{e} conjecture} is that the third term is 1. \\

\noindent{\bf Dedekind's problem.}
The number of monotone Boolean functions of $n$ variables, for $n=1,2,3, \ldots, 8$
(\htmladdnormallink{A7153}{http://www.research.att.com/cgi-bin/access.cgi/as/njas/sequences/eisA.cgi?Anum=007153}):
$$
\begin{array}{l}
1,4,18,166, 7579, 7828352, 2414682040996, \\
56130437228687557907786$$
\end{array}
$$
where the last term was computed by Wiedemann \cite{Wie}.
(As is the case for many of these examples, there are several other versions of this sequence in the database.) \\

\noindent{\bf The Hadamard maximal determinant problem.}
What is the maximal determinant of an $n \times n$ $\{0,1\}$-matrix?
The values for $n=1,2, \ldots, 13$ are (\htmladdnormallink{A3432}{http://www.research.att.com/cgi-bin/access.cgi/as/njas/sequences/eisA.cgi?Anum=003432}):
$$1,1,2,3,5,9, 32, 56, 144, 320, 1458, 3645, 9477,$$
where the last two terms are due to Ehlich, and Ehlich and Zeller \cite{Ehl1},
\cite{EhZ}.
For $n \equiv -1~(\bmod~4)$, Hadamard showed that the $n$th term is equal to 
$$(n+1)^{(n+1)/2} /2^n~,$$
provided that what is now called a Hadamard matrix of order $n$
exists.
In some cases conference matrices give the answer when
$n \equiv 1~(\bmod~4)$, but
the problem of finding the other terms in the sequence seems 
to have been untouched for over 35 years.
It would be nice to have confirmation of the above values as well as
some more terms! \\

\noindent{\bf Enumerating Hadamard matrices.}
The number of Hadamard matrices
of orders $n=4,8,12, \ldots, 28$ (\htmladdnormallink{A7299}{http://www.research.att.com/cgi-bin/access.cgi/as/njas/sequences/eisA.cgi?Anum=007299}) is
$$1,1,1,5,3,60, 487,$$
where the last entry is the work of Kimura \cite{Kimur86},
\cite{Kimur94a},
\cite{Kimur94b},
\cite{KimOh86}.
The {\em Hadamard conjecture} is that such a matrix always exists if $n$ is
a multiple of 4.
Judging by how rapidly these numbers are growing, this should not be hard to prove,
yet it has remained an open question for a century.
Of course, as the example in Section 7 shows, such numerical evidence can be
misleading. \\

\noindent{\bf The kissing number problem.}
How many spheres can touch another sphere of the same size?
For arrangements that occur as part of a lattice packing,
the answers are known for $n=1,2, \ldots, 9$ (\htmladdnormallink{A1116}{http://www.research.att.com/cgi-bin/access.cgi/as/njas/sequences/eisA.cgi?Anum=001116}):
$$2,6,12,24,40, 72, 126, 240 , 272,$$
the last term being due to Watson (see \cite{SPLAG}).
For nonlattice packings, all we know is
$$2,6,12,?,?,?,?,240, \ge 306~.$$
The best bounds known in dimensions 4, 5, 6 and 7 are respectively
$$\mbox{24--25, 40--46, 72--82} \quad\mbox{and}\quad
\mbox{126--140}$$
--- see \cite{SPLAG} for further information. \\

\paragraph{\bf Sorting by prefix reversal.}\label{sorting}
If you can only reverse segments that include the initial term of the current permutation, how many reversals are needed to
transform an arbitrary permutation of $n$ letters to the identity permutation?
To state this another way \cite{AMM75}:
\begin{quote}
The chef in our place is sloppy, and when he prepares a stack of pancakes they come out all different sizes.
Therefore, when I deliver them to a customer, on the way to the table I
rearrange them (so that the smallest winds up on top, and so on, down to the largest at the bottom)
by grabbing several from the top and flipping them over,
repeating this (varying the number I flip) as many times as necessary.
If there are $n$ pancakes, what is the maximum number of flips (as a function
$f(n)$ of $n$) that I will ever have to use to rearrange them?
\end{quote}
The only exact values known are $f(1),\ldots, f(9)$:
$$0,1,3,4,5,7,8,9,10,$$
due to Garey, Johnson and Lin, and Robbins (see \cite{AMM75}, \cite{GaP79}).
It is also known that $f(n) \ge n+1$ for $n \ge 6$, $f(n) 
\ge 17n/16$ if $n$ is a multiple of 16
(so $f(32) \ge 34$), and $f(n) \le (5n+5)/3$, the last two bounds\footnote{(Added May 5, 2000.)
The bound $f(n) \le (5n+5) /3$ was independently obtained by E. Gy\"{o}ri and G. Tur\'{a}n, Stack of pancakes,
{\em Studia Sci. Math. Hungar.}, {\bf 13} (1978), 133--137.
I thank L\'{a}szl\'{o} Lov\'{a}sz for pointing out this reference.}
being due to
Gates and Papadimitriou \cite{GaP79}.
Again it would be nice to have more terms.

Note added Jan. 17, 2001.
John J. Chew III (Department of Mathematics, University of Toronto)
has found that $f(10)$ through $f(13)$ are $11, 13, 14, 15$, respectively,
so the sequence begins:
$$0,1,3,4,5,7,8,9,10,11,13,14,15$$
This has now been added to the database as sequence
\htmladdnormallink{A058986}{http://www.research.att.com/cgi-bin/access.cgi/as/njas/sequences/eisA.cgi?Anum=058986}.

\section{Recursive sequences}
Whereas the sequences in the previous section enumerated some class of objects, the following are self-generated. \\

\noindent{\bf Differences $=$ complement} (\htmladdnormallink{A5228}{http://www.research.att.com/cgi-bin/access.cgi/as/njas/sequences/eisA.cgi?Anum=005228}):
$$1,3,7,12,18,26,35,45,56,69,83,98,114, \ldots$$
The differences (\htmladdnormallink{A30124}{http://www.research.att.com/cgi-bin/access.cgi/as/njas/sequences/eisA.cgi?Anum=030124})
$$2,4,5,6,8,9,10,11,13,14,15,16,17,19, \ldots$$
are the terms not in the sequence!
This is one of many fine self-generating sequences from Hofstadter \cite{GEB}. \\

\noindent{\bf Golomb's sequence.}
The $n$th term is the number of times $n$ appears (\htmladdnormallink{A1462}{http://www.research.att.com/cgi-bin/access.cgi/as/njas/sequences/eisA.cgi?Anum=001462}):
$$1,2,2,3,3,4,4,4,5,5,5,6,6,6,6,7,7,7,7,8, \ldots$$
The $n$th term is the nearest integer to (and converges to)
$$\tau^{2- \tau} n^{\tau -1} ~,$$
where $\tau = (1+ \sqrt{5})/2$ \cite{Gol66},
\cite[Section E25]{UPNT}. \\

\noindent{\bf Wilson's primeth recurrence.}
$a_{n+1}$ is the $a_n$-th prime (\htmladdnormallink{A7097}{http://www.research.att.com/cgi-bin/access.cgi/as/njas/sequences/eisA.cgi?Anum=007097}), shown in Fig.~\ref{FRP1}.
The sequence was sent in by R. G. Wilson V \cite{RGW}, and the last few terms were computed
by P. Zimmermann \cite{Zim1} and M. Del\'{e}glise \cite{Del1}.
Their algorithm is a slightly speeded up version of an algorithm for computing
$\pi (x)$, the number of primes not exceeding $x$, due to
J. C. Lagarias, V. S. Miller and A. M. Odlyzko (see \cite{LO82}).
It is quite remarkable that it is possible to compute so many terms of this
sequence.

\begin{figure}[htb]
$$
\begin{array}{r}
1 \\
2 \\
3 \\
5 \\
11 \\
31 \\
127 \\
709 \\
5381 \\
52711 \\
648391 \\
9737333 \\
174440041 \\
3657500101 \\
88362852307 \\
2428095424619 \\
75063692618249 \\
2586559730396077 \\
\end{array}
$$

\caption{$a_{n+1}$ is the $a_n$-th prime.}
\label{FRP1}
\end{figure}

\noindent{\bf Recam\'{a}n's sequences.}
(i) $a_n = a_{n-1} -n$ if $a_{n-1} -n > 0$ and $a_{n-1} -n$
has not already occurred in the sequence, otherwise $a_n = a_{n-1} +n$ (\htmladdnormallink{A5132}{http://www.research.att.com/cgi-bin/access.cgi/as/njas/sequences/eisA.cgi?Anum=005132}):
$$1,3,6,2,7,13,20,12,21,11,22,10,23,9, \ldots$$
(ii) $a_{n+1} = a_n/n$ if $n$ divides $a_n$, otherwise $a_{n+1} = na_n$ (\htmladdnormallink{A8336}{http://www.research.att.com/cgi-bin/access.cgi/as/njas/sequences/eisA.cgi?Anum=008336}):
$$1,1,2,6,24,120,20,140, 1120, 10080, \ldots$$
These were sent in by B. Recam\'{a}n \cite{Rec1}.
How fast do they grow? \\

\noindent{\bf The \$10,000 sequence.}
In a colloquium talk at AT\&T Bell Labs \cite{Bell}, John Conway
discussed the sequence (\htmladdnormallink{A4001}{http://www.research.att.com/cgi-bin/access.cgi/as/njas/sequences/eisA.cgi?Anum=004001})
$$1,1,2,2,3,4,4,4,5,6,7,7,8,8,8,8,9, \ldots$$
defined by (for $n \ge 3$)
$$a(n+1) = a(a(n)) + a(n+1 - a(n)) ~.$$
(In words, $a(n+1)$ is the $a(n)$th term in
from the left plus the $a(n)$th term in from the right.)
This sequence seems to have been introduced by either David Newman or Douglas Hofstadter around 1986.
In his talk,
Conway said that he could prove that $\frac{a(n)}{n} \to \frac{1}{2}$,
and offered \$10,000 for finding the exact $n$ at which
$\left| \frac{a(n)}{n} - \frac{1}{2} \right|$ last exceeds $\frac{1}{20}$.
My colleague Colin Mallows did not take long to analyze the sequence, and came up with an answer of
6083008742 \cite{CCS1}.

Colin tells me that in fact the problem is actually much easier than either he or John Conway had believed, and the true answer is 1489.
A recent paper by Kubo and Vakil \cite{KV1} also studies this sequence and its generalizations.
See also \cite[Section E31]{UPNT}.

Many variants of this sequence have not yet been analyzed.
Even the rate of growth of this one is not known:
$a(1) = a(2) =1$, $a(n) = a(a(n-2)) + a(n-a(n-2))$, \htmladdnormallink{A5229}{http://www.research.att.com/cgi-bin/access.cgi/as/njas/sequences/eisA.cgi?Anum=005229} \cite{CCS1},
which begins
$$1,1,2,3,3,4,5,6,6,7,7,8,9,10, 10,10, 11,12,12, \ldots$$ \\

\noindent{\bf The Prague clock sequence.}
This is not really recursive
in the same sense as the preceding sequences, but it seems to fit in here.
It is, in any case, delightful.
This is \htmladdnormallink{A28354}{http://www.research.att.com/cgi-bin/access.cgi/as/njas/sequences/eisA.cgi?Anum=028354}:
1, 2, 3, 4, 32, 123, 43, 2123, 432, 1234, 32123, 43212, 34321, 23432, 123432, 
1234321, 2343212, 3432123, 4321234, 32123432, 123432123, 43212343, 
2123432123, 432123432, 1, 2, 3, 4, 32, $\ldots$
According to \cite{Hor}, the sequence indicates how the astronomical clock in Prague strikes the hours.
There is a single bell which
(i) makes from 1 to 4 strokes at a time,
(ii)~the number of strokes follows the sequence
$$\ldots 321 2343 212343 \ldots ~,$$
(iii)~at the $n$th hour, for $n=1,2, \ldots, 24$, the strokes add to $n$, and
(iv)~at the 25th hour there is a single stroke (so the sequence has period 24).
As the reader will see by studying the sequence, its existence depends on two coincidences! \\

\noindent{\bf Cayley's mistake.}
Since the sequences in the database are numbered A1, A2, A3, $\ldots$, several people humorously proposed the ``diagonal'' sequence
(now \htmladdnormallink{A31135}{http://www.research.att.com/cgi-bin/access.cgi/as/njas/sequences/eisA.cgi?Anum=031135}) in which the $n$th term is equal to the $n$th term of $An$:
$$1,2,1,0,2,3,0,6,8,4,63,1,316,42,16, \ldots ,$$
and the even less well-defined sequence (now \htmladdnormallink{A37181}{http://www.research.att.com/cgi-bin/access.cgi/as/njas/sequences/eisA.cgi?Anum=037181})
with $n$th term equal to \htmladdnormallink{A31135}{http://www.research.att.com/cgi-bin/access.cgi/as/njas/sequences/eisA.cgi?Anum=031135}$(n)+1$.
I resisted adding these sequences for a long time, partly out of a desire to maintain
the dignity of the database, and partly because \htmladdnormallink{A22}{http://www.research.att.com/cgi-bin/access.cgi/as/njas/sequences/eisA.cgi?Anum=000022}
was only known to 11 terms!

Sequence \htmladdnormallink{A22}{http://www.research.att.com/cgi-bin/access.cgi/as/njas/sequences/eisA.cgi?Anum=000022} gives the number of ``centered hydrocarbons with $n$ atoms'',
and is based on an 1875 paper of Cayley \cite{Cay75}.
The paper is extremely hard to follow, and gives incorrect values for $n=12$ and 13.
The errors it contains were reproduced in \cite{BuS65}, \cite{HIS} and \cite{EIS},
even though Herrmann \cite{Her80} had pointed out these errors in 1880.
As far as we can tell, a correct verion of this sequence was never published until
Eric Rains and I did so in 1999 \cite{RS99}.
Although Henze and Blair, P\'{o}lya and many others have written articles enumerating related families of chemical compounds (see \cite{RS99} for a brief survey),
this sequence seems to have been forgotten for over 100 years.

Once we determined what it was that Cayley was trying to enumerate,
P\'{o}lya's counting theory quickly gave the answer,
and the correct version of \htmladdnormallink{A22}{http://www.research.att.com/cgi-bin/access.cgi/as/njas/sequences/eisA.cgi?Anum=000022} is now in the database (as
are the two diagonal sequences mentioned above -- mostly because users 
of the database kept proposing them).

\section{Meanders and stamp-folding}
The meandric and stamp-folding numbers are similar to the better-known Catalan numbers (sequence \htmladdnormallink{A108}{http://www.research.att.com/cgi-bin/access.cgi/as/njas/sequences/eisA.cgi?Anum=000108})
in that they are fundamental, easily described combinatorial quantities that arise in
many different parts of mathematics, but differ from them in that there is no
known formula,
and in fact seem to be quite hard to compute.

\begin{figure}[htb]
\centerline{\psfig{file=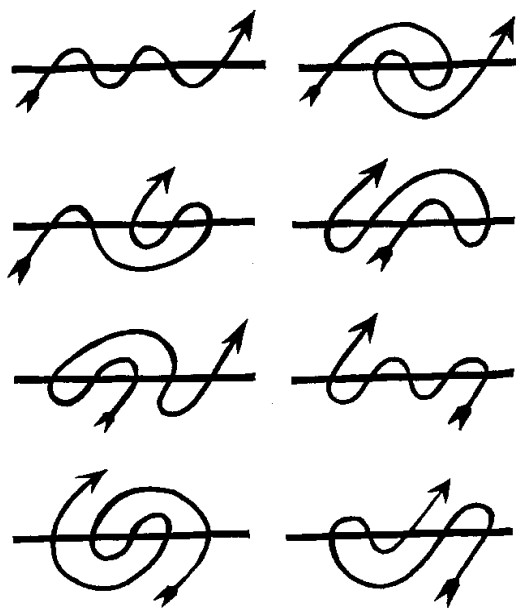,width=3in}}

\caption{The eight ways a river (going from South-West to North-East) can cross a road five times.}
\label{FMe1}
\end{figure}

The stamp-folding problem has a history going back to at least Lucas \cite{Luc}
in the nineteenth century, while the meandric problem
seems to have been first mentioned by Poincar\'{e} \cite{Poin}.
The meandric problem asks:
in how many different ways can a river (starting in the South-West and flowing East)
cross a road $n$ times?
For $n=5$ crossings there are $M_5 =8$ possibilities, shown in Fig. \ref{FMe1}.
The sequence $M_1$, $M_2$, $M_3$, $\ldots$ (\htmladdnormallink{A5316}{http://www.research.att.com/cgi-bin/access.cgi/as/njas/sequences/eisA.cgi?Anum=005316}) begins
$$1,1,2,3,8,14,42,81,262,538,1828, \ldots$$
These are called {\em meandric} numbers, since the river {\em meanders} across the road.
The even-numbered terms $M_2$, $M_4$, $M_6 , \ldots$ give the number of different ways an oriented line can cross a Jordan curve
(\htmladdnormallink{A5315}{http://www.research.att.com/cgi-bin/access.cgi/as/njas/sequences/eisA.cgi?Anum=005315}).
There are several other interpretations, one of which is the number of
``simple alternating transit mazes'' \cite{Ph88}, \cite{Ph89}.

The stamp-folding problem asks the same question, but now the line is only semi-infinite.
Equivalently, how many ways are there to fold a strip of $n$ stamps?
Pictures illustrating the first few terms
of the stamp-folding sequence (\htmladdnormallink{A1011}{http://www.research.att.com/cgi-bin/access.cgi/as/njas/sequences/eisA.cgi?Anum=001011})
can be found on the front cover of \cite{HIS}
and in Figure M4587 of \cite{EIS}.
The sequence begins
$$1,1,2,5,14,38,120, 353,1148, 3527, \ldots$$

No polynomial-time algorithm is known for computing either sequence (on the other hand
it is not known that such algorithms do not exist).
The best algorithms known require on the order of $nC_n$ steps, where $C_n = \frac{1}{n+1} {\binom{2n}{n}}$ is the $n$th
Catalan number (\htmladdnormallink{A108}{http://www.research.att.com/cgi-bin/access.cgi/as/njas/sequences/eisA.cgi?Anum=000108}):
these algorithms are due to Koehler \cite{Koe} for the stamp-folding problem and to
Knuth and Pratt \cite{KP89} and Reeds \cite{RSM} for the meandric problem.
They are too complicated to describe here.

Using these algorithms, Stephane Legendre \cite{Leg} has extended the stamp-folding sequence to 26 terms and the meandric sequence to 25 terms.
Knuth and Pratt \cite{KP89} have computed 17 terms of the $M_2$, $M_4$, $M_6 , \ldots$ subsequence.
Lando and Zvonkin \cite{LZ92}, \cite{LZ93} and Di Francesco, Golinelli and Guitter \cite{FGG96},
\cite{FGG95}, \cite{FGG97} have also studied these sequences.

The exact rate of growth of these sequences is not known.
The best bounds for $M_{2n}$ presently known
appear to be due to Reeds, Shepp and McIlroy \cite{RSM}.
It is easy to see that $M_{2n}$ is submultiplicative, and that $C_n \le M_{2n} \le C_n^2$, which implies that
$$\mu = \lim_{n \to \infty} M_{2n}^{1/n}$$
exists and satisfies $4 \le \mu \le 16$.
In \cite{RSM} it is shown that
$$8.8 \le \mu \le 13.01 ~.$$

Besides the papers already mentioned, the meandric and stamp-folding problems have recently
been discussed by
Arnold \cite{Arn88}, Di~Francesco \cite{DiF96},
Harris \cite{Har98},
Lando and Zvonkin \cite{LZ92}, \cite{LZ93}, Lunnon \cite{Lun68},
Sade \cite{Sade},
Smith \cite{WDS} and
Touchard \cite{Touch50}.

\section{Extremal codes and lattices}
Let $C$ be a binary linear self-dual code of length $n$ in which the weight of every codeword
is a multiple of 4, and let
$$W_C (x,y) = \sum_{c \in C}
x^{n- wt(c)} y^{wt(c)}$$
be its weight enumerator.
Examples are the Hamming code of length 8, with weight enumerator
$$f = x^8 + 14 x^4 y^4 + y^8 ~,$$
and the Golay code of length 24, with weight enumerator
$$g = x^{24} + 759 x^{16} y^8 + 2576 x^{12} y^{12} + 759 x^8 y^{16} + y^{24} ~.$$

A remarkable theorem of Gleason says that the weight enumerator of any 
such code $C$ is a polynomial in $f$ and $g$.
(References for this section are
\cite{SPLAG}, \cite{MS77},
\cite{RS98}.)
E.g. if $C$ has length $n=72$, its weight enumerator can be written as
$$W_C = a_0 f^9 + a_1 f^6 g + a_2 f^3 g^2 + a_3 g^3 ~,$$
for rational numbers $a_0, \ldots, a_3$.
If we choose these numbers so that the minimal distance of
this (hypothetical) code $C$ is as large as possible, we find that
$$W_C =1+ 0x^4 + 0x^8 +0x^{12} +
249849 x^{16} + 18106704 x^{20} + 462962955 x^{24} + \cdots ~,
$$
so that $C$ would have minimal distance 16 (the coefficients form sequence
\htmladdnormallink{A18236}{http://www.research.att.com/cgi-bin/access.cgi/as/njas/sequences/eisA.cgi?Anum=018236}).
The coefficients in this ``extremal'' weight enumerator are all nonnegative integers, but whether a code exists with this weight enumerator is an important unsolved question.

One can perform this calculation for any length that is a multiple of 8, and in \cite{MaS} it was shown that
when as many initial terms as possible are set to zero,
the next term in the extremal weight enumerator is always
positive (so the minimal distance of the hypothetical extremal code is precisely $4[n/24] + 4$).

Lengths that are multiples of 24 are especially interesting.
Figure \ref{FEC1} shows the leading term in the extremal weight enumerator
for a binary self-dual code of length $n=24m$
(this is sequence \htmladdnormallink{A34414}{http://www.research.att.com/cgi-bin/access.cgi/as/njas/sequences/eisA.cgi?Anum=034414}).
\begin{figure}[htb]
$$
\begin{array}{r@{~~~~}r}
n & \multicolumn{1}{c}{\mbox{Coefficient}} \\ \hline
0 & 1 \\
24 & 759 \\
48 & 17296 \\
72 & 249849 \\
96 & 3217056 \\
120 & 39703755 \\
144 & 481008528 \\
168 & 5776211364 \\
192 & 69065734464
\end{array}
$$

\caption{Number of codewords of minimal weight $n/6 +4$ in extremal weight enumerator of length $n$.}
\label{FEC1}
\end{figure}
It is ``obvious'' that these numbers are growing rapidly, and in fact it is shown in \cite{MaS} that if $n=24m$ then the leading
coefficient is
$${\binom{24m}{5}} {\binom{5m-2}{m-1}} / {\binom{4m+4}{5}} ~.$$
The next term in the extremal weight enumerator (the number of codewords
of weight $4[n/24] +8$) also grows rapidly --- see Fig \ref{FEC2} (sequence \htmladdnormallink{A34415}{http://www.research.att.com/cgi-bin/access.cgi/as/njas/sequences/eisA.cgi?Anum=034415}).
\begin{figure}[htb]
$$
\begin{array}{r@{~~~~}r}
n & \multicolumn{1}{c}{\mbox{Coefficient}} \\ \hline
0 & 1 \\
24 & 2576 \\
48 & 535095 \\
72 & 18106704 \\
96 & 369844880 \\
120 & 6101289120 \\
144 & 90184804281 \\
168 & 1251098739072 \\
192 & 16681003659936
\end{array}
$$

\caption{Number of codewords of next-to-minimal weight $n/6+8$ in extremal weight enumerator of
length $n$.}
\label{FEC2}
\end{figure}

Again it is ``obvious'' that these numbers also grow rapidly, especially if one examines a more extensive
table, where the number of  digits continues to increase with 
each term.
Yet --- this shows the danger of drawing conclusions just from numerical data ---
it is proved in \cite{MOS} that this sequence goes negative at $n= 3696$ and stays
negative forever.
(See Theorem 29 of \cite{RS98} for more precise information.)

So codes corresponding to these extremal weight enumerators
certainly do not exist for $n \ge 3696$ (since the
weight enumerator of a genuine code must have nonnegative
coefficients).
They are known to exist for $n=24$ (the Golay code)
and 48 (a quadratic residue code),
but every case from 72 to 3672 is open.
For further details, and in particular for a description
of the analogous situation when $n$ is a multiple of 8 but not of 24, see \cite{RS98}.

The situation for lattices is similar.
Consider an even unimodular lattice $\Lambda$ of dimension $n$, with theta series
$$\Theta_\Lambda (q) = \sum_{u \in \Lambda} q^{u \cdot u} ~.$$
Examples are the eight-dimensional root lattice $E_8$, with theta series
\begin{eqnarray*}
f & = & 1 + 240 q^2 + 2160q^4 + 6720 q^6 + \cdots \\
& = & 1+ 240 \sum_{m=1}^\infty \sigma_3 (m) q^{2m} ~,
\end{eqnarray*}
where $\sigma_3 (m)$ is the sum of the cubes of the divisors of $m$,
and the
24-dimensional Leech lattice, with theta series
\begin{eqnarray*}
g & = & 1+ 196560 q^4 + 16773120q^6 +
398034000 q^8 + \cdots \\
& = & f^3 - 720 q^2 \prod_{m=1}^\infty (1-q^{2m} )^{24} ~.
\end{eqnarray*}
The Taylor series expansion of the last product defines the famous Ramanujan numbers.
(The coefficients of $f$ and $g$ give sequences \htmladdnormallink{A4009}{http://www.research.att.com/cgi-bin/access.cgi/as/njas/sequences/eisA.cgi?Anum=004009} and \htmladdnormallink{A8408}{http://www.research.att.com/cgi-bin/access.cgi/as/njas/sequences/eisA.cgi?Anum=008408}, respectively;
the Ramanujan numbers form sequence \htmladdnormallink{A594}{http://www.research.att.com/cgi-bin/access.cgi/as/njas/sequences/eisA.cgi?Anum=000594}.)

A theorem due essentially to Hecke says that the theta series of any such lattice
$\Lambda$ is a polynomial in $f$ and $g$.
We may then define extremal theta series just as we defined extremal weight enumerators.
It is known \cite{MOS} that the leading term
in the extremal theta series is positive, but that
again the next coefficient eventually becomes negative (even
though the sequence begins with about 1700 exponentially
growing terms).
The dimension where the first negative coefficient
appears is around 41000.
The leading terms are shown in Fig.~\ref{FEC3} (sequence \htmladdnormallink{A34597}{http://www.research.att.com/cgi-bin/access.cgi/as/njas/sequences/eisA.cgi?Anum=034597}).
\begin{figure}[htb]
$$
\begin{array}{r@{~~~~}r}
n & \multicolumn{1}{c}{\mbox{Coefficient}} \\ \hline
0 & 1 \\
24 & 196560 \\
48 & 52416000 \\
72 & 6218175600 \\
96 & 565866362880 \\
120 & 45792819072000
\end{array}
$$

\caption{Number of shortest vectors $u$ (with $u \cdot u = n/12 +2$) in extremal theta series at dimension $n$.}
\label{FEC3}
\end{figure}

Lattices corresponding to the extremal theta series exist
for $n=24$ (the Leech lattice) and 48 (at least three
inequivalent lattices exist), but every
case from 72 to about 41000 is open.
The 24- and 48-dimensional extremal lattices
provide record kissing numbers in those dimensions;
it would be nice to have a 72-dimensional example.

\section{Levine's sequence}
In the summer of 1997 Lionel Levine \cite{Lev97}
submitted a new sequence to the table,
a sequence of such beauty that it immediately caught
the attention of several of my colleagues.
It is constructed via the array in Fig.~\ref{FL}.
\begin{figure}[htb]
$$
\begin{array}{|ccccccccccccccccccccccc|} \hline
&&&&&&&&&&&&&&&&&&&&&1 & \underline{1} \\ \hline
&&&&&&&&&&&&&&&&&&&&&1 & \underline{2} \\ \hline
&&&&&&&&&&&&&&&&&&&&1 & 1 & \underline{2} \\ \hline
&&&&&&&&&&&&&&&&&&&1 & 1 & 2 & \underline{3} \\ \hline
&&&&&&&&&&&&&&&&1 & 1 & 1 & 2 & 2 & 3 & \underline{4} \\ \hline
&&&&&&&&&1 & 1 & 1 & 1 & 2 & 2 & 2 & 3 & 3 & 4 & 4 & 5 & 6 & \underline{7} \\ \hline
1& 1 & 1 & 1 & 1 & 1 & 1 & 2 & 2 & 2 & 2 & 2 & 2 & 3 & 3 & 3 & 3 & 3 & 4 & 4 & ~ & ~ & ~ \\
&& 4 & 4 & 5 & 5 & 5 & 5 & 6 & 6 & 6 & 7 & 7 & 7 & 8 & 8 & 9 & 9 & 10 & 11 & 12 & 13 & \underline{14} \\  \hline
\end{array}
$$
\caption{The array that produces Levine's sequence
$1,2,2,3,4,7,14, \ldots$.}
\label{FL}
\end{figure}
If a row of the array contains the numbers
$$a_1 ~~a_2 ~~ a_3 ~ \ldots ~ a_{k-1} ~~ a_k$$
then the next row contains
$$
\mbox{$a_k$  1's,}
\quad
\mbox{$a_{k-1}$  2's,} \quad
\mbox{$a_{k-2}$  3's, $\ldots$}
$$
Levine's sequence (\htmladdnormallink{A11784}{http://www.research.att.com/cgi-bin/access.cgi/as/njas/sequences/eisA.cgi?Anum=011784}) is obtained by taking the last term in each row:
\begin{equation}\label{EqLe1}
\begin{array}{l}
1, 2, 2, 3, 4, 7, 14, 42, 213, 2837, 175450, 
139759600, 6837625106787, \\
266437144916648607844, 
508009471379488821444261986503540, ...
\end{array}
\end{equation}

The terms grow unexpectedly rapidly!
The $n$th term $L_n$ is
\begin{itemize}
\item[(i)]
the sum of the elements in row $n-2$
\item[(ii)]
the number of elements in row $n-1$
\item[(iii)]
the last element in row $n$
\item[(iv)]
the number of 1's in row $n+1$

$\ldots$
\end{itemize}
Furthermore, if $s(n,i)$ denotes the sum of the first $i$ elements in
row $n$, then we have
\begin{itemize}
\item[(v)]
$L_{n+2} = s (n, L_{n+1} )$
\item[(vi)]
$L_{n+3} = \displaystyle\sum_{i=1}^{L_{n+1}} s(n,i)$
\item[(vii)]
$L_{n+4} = \displaystyle\sum_{i=1}^{L_{n+1}} {\displaystyle\binom{s(n,i)+1}{2}}$.
\end{itemize}
The latter identity was found by Allan Wilks
\cite{Wilks}, who also found a more complicated formula for
$L_{n+5}$, and used it to compute the last two terms
shown in \eqn{EqLe1}.
No other terms are known!

As to the rate of growth, we have only a crude estimate.
Bjorn Poonen and Eric Rains \cite{PoRa97} showed that
\begin{equation}\label{EqL3}
\log L_n \sim c \tau^n ~,
\end{equation}
where $\tau = (1+ \sqrt{5})/2$.
{\em Sketch of Proof.}
(a)~$L_{n+2} \le L_{n+1} L_n$, and so $\log L_n$ is bounded above by a Fibonacci-like sequence.
(b)~The sum of the $(n+1)$st row is at most
$$
{\binom{\left[ \frac{L_{n+2}}{L_n} \right] +1}{2}} L_n ~,
$$
which implies
$$\frac{L_{n+3}}{2L_{n+2}} \ge \frac{L_{n+2}}{2L_{n+1}} \frac{L_{n+1}}{2L_n}$$
and so $\log (L_{n+1} / 2L_n )$ is bounded below by a Fibonacci-like sequence.

Colin Mallows \cite{clm} has determined numerically that a reasonably
good approximation to $L_n$ is given by
$$\frac{1}{c_1} e^{c_2 \tau^n}$$
where $c_1 \approx 0.277$, $c_2 \approx 0.05427$.
It would be nice to have better estimates for $L_n$, and one or more
additional terms.

\section{Arrays of Numbers}
Besides number {\em sequences}, the database also contains {\em arrays}
of numbers that have been converted to sequences.
Triangular arrays are read by rows, in the obvious way.
E.g. Pascal's triangle of binomial coefficients
$$
\begin{array}{ccccccccc}
~ & ~ & ~ & ~ & 1 \\
~ & ~ & ~ & 1 & ~ & 1 \\
~ & ~ & 1 & ~ & 2 & ~ & 1 \\
~ & 1 & ~ & 3 & ~ & 3 & ~ & 1 \\
1 & ~ & 4 & ~ & 6 & ~ & 4 & ~ & 1 \\
\end{array}
$$
becomes sequence \htmladdnormallink{A7318}{http://www.research.att.com/cgi-bin/access.cgi/as/njas/sequences/eisA.cgi?Anum=007318}:
$$1,1,1,1,2,1,1,3,3,1,1,4,6,4,1, \ldots$$

Square arrays are read by antidiagonals, usually in this order:
$$
\begin{array}{lllll}
a_0 & a_2 & a_5 & a_9 & \ldots \\
a_1 & a_4 & a_8 & \multicolumn{2}{l}{\ldots} \\
a_3 & a_7 & \multicolumn{3}{l}{\ldots} \\
a_6 & \multicolumn{4}{l}{\ldots} \\
\multicolumn{5}{l}{\ldots}
\end{array}
$$
E.g. the Nim-addition table \cite{ONAG}
$$\begin{array}{cccccc}
0 & 1 & 2 & 3 & 4 & \ldots \\
1 & 0 & 3 & 2 & 5 & \ldots \\
2 & 3 & 0 & 1 & 6 & \ldots \\
3 & 2 & 1 & 0 & 7 & \ldots \\
\cdot & \cdot & \cdot
\end{array}
$$
becomes sequence \htmladdnormallink{A3987}{http://www.research.att.com/cgi-bin/access.cgi/as/njas/sequences/eisA.cgi?Anum=003987}:
$$0,1,1,2,0,2,3,3,3,3,4,2,0,2,4, \ldots$$

Other classical arrays are the Stirling numbers of both kinds, Eulerian numbers,
etc.

A less well-known array arises from {\bf Gilbreath's conjecture}.
This conjecture states that if one writes down the primes in a row, and underneath the absolute values of the differences, as in Fig. \ref{FG},
then the leading terms (shown $\underline{\mbox{underlined}}$)
of all rows except the first are equal to 1 (\cite{UPNT}, \S A10).
The corresponding sequence (\htmladdnormallink{A36262}{http://www.research.att.com/cgi-bin/access.cgi/as/njas/sequences/eisA.cgi?Anum=036262}) is
$$2,1,3,1,2,5,1,0,2,7,1,2,2,4,11,1, \ldots$$
Odlyzko \cite{Od93} has verified the conjecture out to $3 \times 10^{11}$.
\begin{figure}[htb]
$$
\begin{array}{ccccccccccccccccccc}
2~~~ & ~~~~ & 3~~~  & ~~~~ & 5~~~ & ~ & 7~~~ & ~~~~ & 11~~ & ~ & 13~~ & ~ & 17~~ & ~ & 19~~ & ~ & 23~~ & ~~~~ & \ldots \\
~ & \underline{1} & ~ & 2 & ~ & 2 & ~ & 4 & ~ & 2 & ~ & 4 & ~ & 2 & ~ & 4 & ~ & \ldots \\
~ & ~ & \underline{1} & ~ & 0 & ~ & 2 & ~ & 2 & ~ & 2 & ~ & 2 & ~ & 2 & ~ & \ldots \\
~ & ~ & ~ & \underline{1} & ~ & 2 & ~ & 0 & ~ & 0 & ~ & 0 & ~ &0 & ~ & \ldots \\
~ & ~ & ~ & ~ & \underline{1} & ~ & 2 & ~ & 0 & ~ & 0 & ~ & 0 & ~ & \ldots \\
~ & ~ & ~ & ~ & ~ &  \underline{1} & ~ & 2 & ~ & 0 & ~ & 0 & ~ & \ldots \\
\multicolumn{16}{c}{\ldots}
\end{array}
$$
\caption{Gilbraith's conjecture is that the leading
terms of all rows in this array
except the first are always 1 (the top row contains the primes, subsequent rows are the absolute
values of the differences of the previous row).}
\label{FG}
\end{figure}

\section{The Wythoff array}
This array shown in Fig. \ref{FW} has many wonderful properties, some of which are mentioned here.
I learned about most of these properties from John Conway
\cite{JHC}, but this array has a long history
--- see Fraenkel and Kimberling \cite{FrKi94},
Kimberling \cite{Kim93},
\cite{Kim93a}, \cite{Kim94},
\cite{Kim95}, \cite{Kim95a},
\cite{Kim95b},
\cite{Kim97}, Morrison \cite{Mor80} and Stolarsky \cite{Sto76},
\cite{Sto77}.
It is related to a large number of sequences in the database (the main
entry is \htmladdnormallink{A35513}{http://www.research.att.com/cgi-bin/access.cgi/as/njas/sequences/eisA.cgi?Anum=035513}).
\begin{figure}[htb]
$$
\begin{array}{rr|rrrrrrrrr}
0 & 1 & 1 & 2 & 3 & 5 & 8 & 13 & 21 & 34 & 55 \\
1 & 3 & 4 & 7 & 11 & 18 & 29 & 47 & 76 & .. \\
2 & 4 & 6 & 10 & 16 & 26 & 42 & 68 & .. \\
3 & 6 & 9 & 15 & 24 & 39 & 63 & .. \\
4 & 8 & 12 & 20 & 32 & 52 & .. \\
5 & 9 & 14 & 23 & 37 & 60 & .. \\
6 & 11 & 17 & 28 & 45 & 73 & .. \\
7 & 12 & 19 & 31 & 50 & 81 & .. \\
. & .. & .. & .. & .. & .. & ..
\end{array}
$$
\caption{The Wythoff array.}
\label{FW}
\end{figure}

\paragraph{Construction (1).}
The two columns to the left of the vertical
line consist respectively of the nonnegative
integers $n$, and the
{\em lower Wythoff sequence}
(\htmladdnormallink{A201}{http://www.research.att.com/cgi-bin/access.cgi/as/njas/sequences/eisA.cgi?Anum=000201}), whose $n$th term is $[(n+1) \tau ]$.
The rows are then filled in by the Fibonacci rule that each term is the sum of the two previous terms.
The entry $n$ in the first column is the {\em index} of that row.

\paragraph{Definition.}
The {\em Zeckendorf expansion} of a number $n$ is obtained by repeatedly subtracting the
largest possible Fibonacci number until nothing
remains.
E.g. $100 = 89+8+3= F_{11} + F_6 + F_4$.
The {\em Fibonacci successor} to $n$, $Sn$, say, is found by replacing
each $F_i$ in the Zeckendorf expansion by $F_{i+1}$.
E.g. the Fibonacci successor to 100 is $S100 = F_{12} + F_7 + F_5 = 144 + 13 + 5 = 162$.

\paragraph{Construction (2).}
The two columns to the left of the vertical line in Fig. \ref{FW} read $n$, $1+Sn$;
then after the vertical line the row continues
$$
m \quad Sm \quad SSm \quad SSSm \quad SSSSm \quad \ldots ~,$$
where $m=n+1 +Sn$.
\paragraph{Construction (3).}
Let $\{S1, S2, S3, \ldots \} = \{2,3,5,7,8,10,11, \ldots \}$ be the sequence of
Fibonacci successors (\htmladdnormallink{A22342}{http://www.research.att.com/cgi-bin/access.cgi/as/njas/sequences/eisA.cgi?Anum=022342}).
The first column to the right of the line consists of the numbers not in that sequence: 1, 4, 6, 9, 12, $\ldots$ (\htmladdnormallink{A7067}{http://www.research.att.com/cgi-bin/access.cgi/as/njas/sequences/eisA.cgi?Anum=007067}).
The rest of each row is filled in by repeatedly applying $S$.
\paragraph{Construction (4).}
The entry in row $n$ and column $k$ is
$$[(n+1) \tau F_{k+2} ] + F_{k+1} n$$
(where $k=0$ indicates the first column to the right of the vertical line).

Some properties of the array to the right of the line
are the following:
\begin{itemize}
\item[(i)]
The first row consists of the
Fibonacci sequence $1,2,3,5,8, \ldots$
\item[(ii)]
Every row satisfies the Fibonacci recurrence.
\item[(iii)]
The leading term in each row is the smallest number not found
in any earlier row.
\item[(iv)]
Every positive integer appears exactly once.
\item[(v)]
The terms in any row or column are monotonically increasing.
\item[(vi)]
Every positive Fibonacci-type sequence (i.e. satisfying
$a(n)=a(n-1)+a(n-2)$ and eventually positive) appears
as some row of the array.
\item[(vii)]
The terms in any two adjacent rows alternate.
\end{itemize}
There are infinitely many arrays with Properties 1--7, see \cite{Kim95a}.

The $n$th term of the {\em vertical para-Fibonacci sequence}
$$
0, 0, 0, 1, 0, 2, 1, 0, 3, 2, 1, 4, 0, 5, 3, 2, 6, 1, 7, 4, 0, 8, 5,  \ldots
$$
(\htmladdnormallink{A19586}{http://www.research.att.com/cgi-bin/access.cgi/as/njas/sequences/eisA.cgi?Anum=019586})
gives the index (or parameter) of the row of the Wythoff array that contains $n$.
This sequence also has some nice fractal-like properties:

(a)~If you delete the first occurrence of each number, the
sequence is unchanged.
Thus if we delete the underlined numbers from
$$
\underline{0}, 0, 0, \underline{1}, 0, \underline{2}, 1, 0,
\underline{3}, 2, 1, \underline{4}, 0, \underline{5},
3, 2, \underline{6}, 1, \underline{7}, 4, 0, \underline{8}, 5, \ldots
$$
we get
$$
0, 0, 0, 1, 0, 2, 1, 0, 3, 2, 1, 4, 0, 5, 3, 2, 6, 1, 7, 4, 0, 8, 5, \ldots
$$
again!

(b)~Between any two consecutive  0's
we see a permutation of the first few positive integers, and these
nest, so the sequence can be rewritten as (read across the rows):
$$
\begin{array}{ccccccccccccc}
0 \\
0 \\
0& ~ & ~ & ~ & ~ & ~ & ~ && 1 \\
0& ~ & ~ & ~ & ~ & 2 & ~ & ~ & 1 \\
0& ~ & ~ & 3 & ~ & 2 & ~ & ~ & 1 & ~ & ~ & 4 \\
0& ~ & 5 & 3 & ~ & 2 & ~ & 6 & 1 & ~ & 7 & 4 \\
0 & 8 & 5 & 3 & 9 & 2 & 10 & 6 & 1 & 11 & 7 & 4 & 12
\end{array}$$

The $n$th term of the {\em horizontal para-Fibonacci sequence}
$$
1, 2, 3, 1, 4, 1, 2, 5, 1, 2, 3, 1, 6, 1, 2, 3, 1, 4, 1, 2, 7, 1, 2,  \ldots
$$
(\htmladdnormallink{A35612}{http://www.research.att.com/cgi-bin/access.cgi/as/njas/sequences/eisA.cgi?Anum=035612})
gives the index (or parameter) of the column of the Wythoff
array that contains $n$.
This sequence also has some nice properties.

I hope I have said enough to convince you that the Wythoff array is
well worth studying and full of surprises.

\section{The Boustrophedon transform}
The Taylor series for $\sin x$ and $\cos x$ are easily remembered,
but most people have trouble with
\begin{eqnarray*}
\tan x & = &
1x + 2 \frac{x^3}{3!} + 16 \frac{x^5}{5!} + 272 \frac{x^7}{7!} + \cdots ~, \\
{\rm sec} ~x & = & 1 + 1 \frac{x^2}{2!} +
5 \frac{x^4}{4!} + 61 \frac{x^6}{6!} + \cdots ~.
\end{eqnarray*}
However, their coefficients can be calculated from the array in Fig. \ref{FB}.
\begin{figure}[htb]
$$
\begin{array}{ccccccccccccccc}
~ & ~ & ~ & ~ & ~ & ~ & ~ & 1 \\
~ & ~ & ~ & ~ & ~ & ~ & 0 & ~ & 1 \\
~ & ~ & ~ & ~ & ~ & 1 & ~ & 1 & ~ & 0 \\
~ & ~ & ~ & ~ & 0 & ~ & 1 & ~ & 2 & ~ & 2 \\
~ & ~ & ~ & 5 & ~ & 5 & ~ & 4 & ~ & 2 & ~ & 0 \\
~ & ~ & 0 & ~ & 5 & ~ & 10 & ~ & 14 & ~ & 16 & ~ & 16 \\
~ & 61 & ~ & 61 & ~ & 56 & ~ & 46 & ~ & 32 & ~ & 16 & ~ & 0 \\
0 & ~ & 61 & ~ & 122 & ~ & 178 & ~ & 224 & ~ & 256 & ~ & 272 & ~ & 272 \\
\multicolumn{15}{c}{\cdots}
\end{array}
$$
\caption{The secant-tangent triangle.}
\label{FB}
\end{figure}
The nonzero entries on the left are the {\em secant numbers} (\htmladdnormallink{A364}{http://www.research.att.com/cgi-bin/access.cgi/as/njas/sequences/eisA.cgi?Anum=000364}):
$$1,1,5,61,1385,50521, 2702765, \ldots ,$$
and those on the right are the
{\em tangent numbers} (\htmladdnormallink{A182}{http://www.research.att.com/cgi-bin/access.cgi/as/njas/sequences/eisA.cgi?Anum=000182})
$$1,2, 16, 272, 7936, 353792, 22368256, \ldots$$
while the combination of the two sequences (\htmladdnormallink{A111}{http://www.research.att.com/cgi-bin/access.cgi/as/njas/sequences/eisA.cgi?Anum=000111}):
\begin{equation}\label{EqB1}
1,1,1,2,5,16,61,272,1385, 7936, \ldots
\end{equation}
are usually called the {\em Entringer numbers}.
The latter count permutations of
\linebreak
$\{1,2, \ldots, n \}$ that
alternately fall and rise.

This array is filled by a rule somewhat similar to that for Pascal's
triangle:
the rows are scanned alternately from right to left
and left to right,
the leading entry in each row is 0,
and every subsequent entry is the sum of the previous
entry in the same row and the entry above it in the previous row.
(This is the {\em boustrophedon} or ``ox-plowing'' rule.)
The earliest reference I have seen to this triangle is Arnold \cite{Arn91},
who calls it the Euler-Bernoulli
triangle.
However, it may well be much older origin.
\cite{bous} gives many other references.

Richard Guy \cite{Guy95} observed that if the entries at the beginnings of the rows
are changed from $1,0,0, \ldots$ to
say $1,1,1,1,1, \ldots$ or to $1,2,4,8,16, \ldots$ then the
numbers that appear at the ends of the rows form interesting-looking
sequences that were not to be found in \cite{EIS}, and asked if they had a combinatorial interpretation.
Using $1,1,1, \ldots$ for example the triangle becomes
$$
\begin{array}{ccccccccccc}
~ & ~ & ~ & ~ & ~ & 1 \\
~ & ~ & ~ & ~ & 1 & ~ & 2 \\
~ & ~ & ~ & 4 & ~ & 3 & ~ & 1 \\
~ & ~ & 1 & ~ & 5 & ~ & 8 & ~ & 9 \\
~ & 24 & ~ & 23 & ~ & 18 & ~ & 10 & ~ & 1 \\
1 & ~ & 25 & ~ & 48 & ~ & 66 & ~ & 76 & ~ & 77 \\
\multicolumn{11}{c}{\cdots}
\end{array}
$$
yielding the sequence (\htmladdnormallink{A667}{http://www.research.att.com/cgi-bin/access.cgi/as/njas/sequences/eisA.cgi?Anum=000667})
\begin{equation}\label{EqB2}
1,2,4,9,24,77, 294, 1309, \ldots
\end{equation}
We may regard this process as carrying out a transformation (the
{\em boustrophedon transform})
of sequences: if
the numbers at the beginnings of the rows are $a_0$, $a_1$, $a_2 , \ldots$
(the input sequence) then the numbers at the ends of the rows,
$b_0$, $b_1$, $b_2, \ldots$ (say) are the output
sequence.
In \cite{bous} we showed that there is a simple relationship
between the input and
output sequences:
their exponential generating functions
$${\cal A} (x) = \sum_{n=0}^\infty a_n \frac{x^n}{n!} , \quad
{\cal B} (x) = \sum_{n=0}^\infty b_n \frac{x^n}{n!}
$$
are related by
$${\cal B} (x) = ({\rm sec}~x + \tan x) {\cal A} (x) ~.$$
We also give a combinatorial interpretation of the $\{b_n\}$.
E.g. in \eqn{EqB2}, $b_n$ is the number of up-down subsequences
of $\{1, \ldots, n\}$, so that $b_3 =9$ corresponds to
$\emptyset ,1,2,3,12,13,23,132,231$.

The Entringer sequence \eqn{EqB1} then has the property that it shifts
one place left under the boustrophedon transform.
The lexicographically earliest sequence that shifts {\em two}
places left under this transform (\htmladdnormallink{A661}{http://www.research.att.com/cgi-bin/access.cgi/as/njas/sequences/eisA.cgi?Anum=000661}) is
$$1,0,1,1,2,6,17, 62, 259, 1230, 6592, \ldots$$
We do not know what this enumerates!

Many examples of similar ``eigen-sequences'' for other transformations
of sequences can be found in Donaghey \cite{MB13},
Cameron \cite{MB6}, and especially \cite{eigen}.
E.g. the sequence giving the number of planted
achiral trees \cite{MB19}, \cite{MB25} (\htmladdnormallink{A3238}{http://www.research.att.com/cgi-bin/access.cgi/as/njas/sequences/eisA.cgi?Anum=003238}):
$$1,1,2,3,5,6,10,11,16,19,26, \ldots$$
has the property that it shifts left one place under the ``inverse M\"{o}bius transformation'' given by
$$b_n = \sum_{d|n} a_d ~.$$

\section{Tchoukaillon solitaire (or Mancala, or Kalahari)}
\begin{figure}[htb]
\centerline{\psfig{file=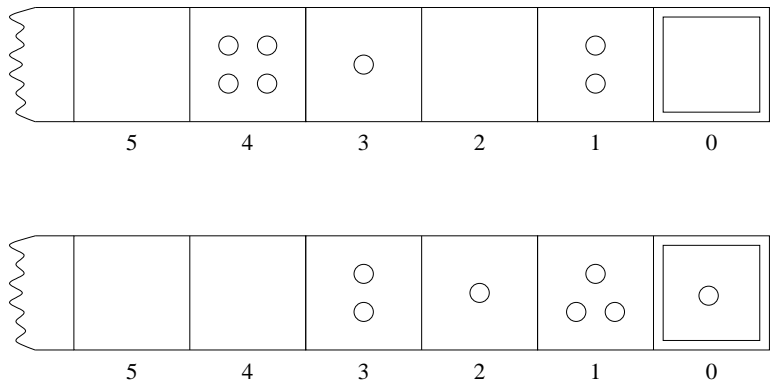,width=4in}}
\caption{A move in Tchoukaillon solitaire.}
\label{FTc1}
\end{figure}

These are ancient board games, with hundreds of variants
and many different names.
The version to be described here is called Tchoukaillon solitaire.
It has been studied by several authors (see for example
Betten \cite{Bet88} and Broline and Loeb \cite{BL95}).
It is played on a board with a row of holes numbered
$0,1,2, \ldots$ (see Fig. \ref{FTc1}).

The game begins with $n$ stones placed anywhere except in hole 0.
A move consists in picking up the stones in some hole and placing one in each lower-numbered hole.
If the last stone falls in hole 0 then play continues, otherwise the game is lost.
The objective is to get all the stones into hole 0.

The game is interesting because there is a unique winning position for any number of stones.
These winning positions are shown in Fig. \ref{FTc2} (sequence \htmladdnormallink{A28932}{http://www.research.att.com/cgi-bin/access.cgi/as/njas/sequences/eisA.cgi?Anum=028932}), and can be found by playing the game backwards.
\begin{figure}[htb]
$$
\begin{array}{r@{~~~~}r}
n & \multicolumn{1}{c}{\mbox{Position}} \\ \hline
0 & 0 \\
1 & 1 \\
2 & 20 \\
3 & 21 \\
4 & 310 \\
5 & 311 \\
6 & 4200 \\
7 &4201 \\
8 &4220 \\
9 &4221 \\
10 &53110 \\
11 &53111 \\
12 &642000 \\
13 & 642001 \\
.. & ...~~~~
\end{array}
$$

\caption{The unique winning position for $n$ stones in Tchoukaillon solitaire.}
\label{FTc2}
\end{figure}

The array can be more explicitly constructed by the rule that if the first 0
in a row (counting from the right)
is in position $i$, then the next row is obtained by writing $i$ in position $i$ and subtracting 1 from all earlier positions.
The sequence of successive values of $i$ (\htmladdnormallink{A28920}{http://www.research.att.com/cgi-bin/access.cgi/as/njas/sequences/eisA.cgi?Anum=028920}) is
$$1,2,1,3,1,4,1,2,1,5,1,6,1,2, \ldots$$

Let $t(k)$ denote the position where $k$ occurs for the first time in this sequence.
The values of $t(1)$, $t(2)$, $t(3), \ldots$ are (sequence A2491):
$$1,2,4,6,10,12,18,22,30,34,42, \ldots$$
This sequence has some very nice properties.
It has been investigated by (in addition to the references mentioned above)
David \cite{Dav57}, Erd\H{o}s and Jabotinsky \cite{ErJa58} and
Smarandache \cite{smar}.

(i)
$t(n)$ can be obtained by starting with $n$ and successively
rounding up to the next multiple of $n-1$, $n-2, \ldots, 2,1$.
E.g. if $n=10$, we obtain
$$10 \to 18 \to 24 \to 28 \to 30 \to 30 \to 32 \to 33 \to 34 \to 34 ~,$$
so $t(10) =34$.

\begin{figure}[htb]
\centerline{\psfig{file=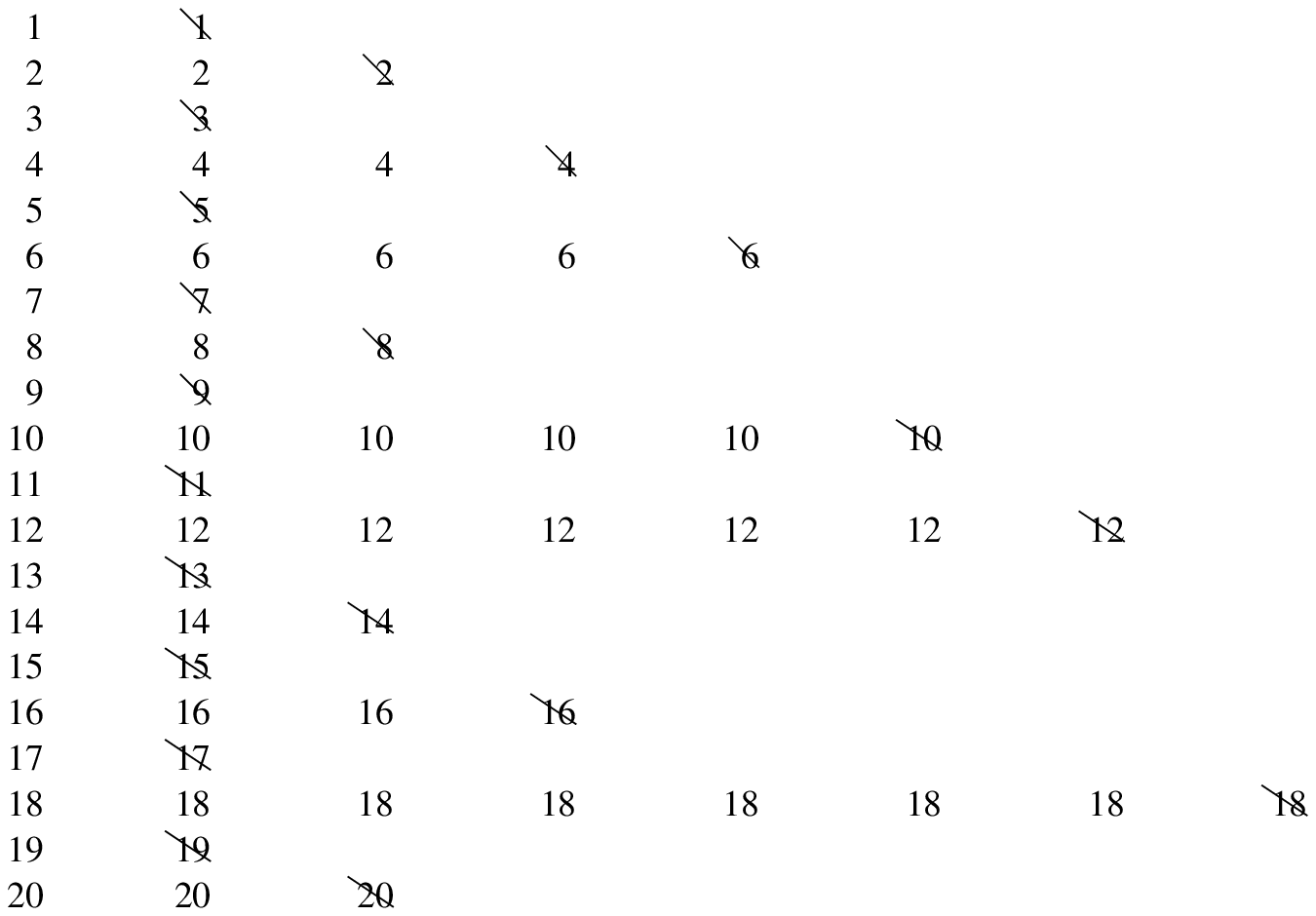,width=4in}}
\caption{A sieve to generate the sequence $t(1), t(2), \ldots = 1,2,4,6,10,12,18, \ldots$.
At stage $n$, the first number and every $n$th are crossed off.}
\label{FTc3}
\end{figure}

(ii)
The sequence can be obtained by a sieving process:
write $1,2, \ldots$ in a column.
To get the second column, cross off 1 and every second number.
To get the third column, cross off the first and every third number.
Then cross off the first and every fourth number, and so on (see Fig. \ref{FTc3}).
The top number in column $n$ is $t(n)$.
Comparison of Figures \ref{FTc2} and \ref{FTc3} shows that connection with the solitaire game.

(iii)
Finally, Broline and Loeb \cite{BL95} (extending the work of
the other authors mentioned) show that, for large $n$,
$$t(n) = \frac{n^2}{\pi} + O(n) ~.$$
It is a pleasant surprise to see $\pi$ emerge from such a simple game.

\end{document}